\documentclass[preprint,12pt]{elsarticle}

\usepackage{amssymb}
\usepackage{amsmath}
\newtheorem{thm}{Theorem}

\newdefinition{rmk}{Remark}
\newdefinition{defn}{Definition}
\newproof{pf}{Proof}
\newproof{pr}{Proof}
\newproof{pot}{Proof of Theorem \ref{thm2}}
\newtheorem{alg}[thm]{Algorithm}
\newcommand{\norm}[1]{\left\lVert#1\right\rVert}

\journal{}

\begin{document}

\begin{frontmatter}



%
 \title{An   upper $J$- Hessenberg reduction of a matrix through symplectic Householder  transformations\\
}
 \author[AS]{Ahmed Salam\corref{cor1}}
 \ead{ahmed.salam@univ-littoral.fr}
 \author[AS,HB]{Haithem Ben Kahla}
 \ead{benkahla@univ-littoral.fr} 
 \cortext[cor1]{corresponding author}

 \address[AS]{University Lille Nord de France, ULCO, LMPA. BP 699, F62228, Calais, Cedx, France.}
 \address[HB]{University of Tunis El Manar, ENIT-LAMSIN, BP 37, 1002, Tunis, Tunisia.}



\begin{abstract}
In this paper, we introduce a reduction of a matrix  to a condensed form, the upper $J$- Hessenberg form, via  elementary symplectic Householder  transformations, which are rank-one modification of the identity . Features of the reduction are highlighted.  Two variants numerically more stables are then derived. 
Some numerical experiments  are given,  showing the  efficiency of these variants.
\end{abstract}


\end{frontmatter}


\section{Introduction}

Let $A$ be a $2n\times 2n$ real matrix. The $SR$ factorization 
consists in writing $A$ as a product $SR$, where $S$ is
symplectic and $R=\left[
\begin{array}{ll}
R_{11} & R_{12}\\
 R_{21} & R_{22}
\end{array}
 \right]$
 is such that $R_{11},\;R_{12},\;R_{22}$ are upper triangular and
 $R_{21}$ is strictly upper triangular
  \cite{ Buns1, Dela2}. This decomposition  plays an important role
 in structure-preserving methods for solving the
eigenproblem of a class of structured matrices.

More precisely, the $SR$ decomposition  can be interpreted  as the analog of the $QR$
decomposition \cite{Golu}, when instead of an Euclidean space, one considers a symplectic space : a 
linear space, equipped with a skew-symmetric inner product (see
for example \cite{Sal1} and the references therein). The 
 orthogonal group with respect to this indefinite inner product,  is  called the 
symplectic group and is unbounded  (contrasting with the Euclidean case).


There are two classes of methods for computing the  $SR$ decomposition.  The first lies in the Gram-Schmidt like algorithms and leads 
 to the 
 symplectic Gram-Schmidt (SGS) algorithms. The second class  is constructed from  a variety of elementary symplectic transformations. Each choice of such transformations leads to the corresponding $SR$ decomposition. Since these elementary transformations  are  quite heterogeneous,  the $SR$ decomposition is considerably  affected by their choice. 

   Results
  on numerical aspects of SGS-algorithms can be found
  for example in   \cite{Sal1}. These
   algorithms and their modified versions are usually involved in
structure-preserving Krylov subspace-type methods, for
 sparse and large structured matrices.


In the literature, the symplectic elementary transformations involved in the $SR$ decomposition can be partitioned in two subsets. The first subset is constituted  of two kind of both symplectic and orthogonal transformations introduced in  \cite{Paig, Vloa}
  and a third symplectic
    but non-orthogonal  transformations,  proposed in
     \cite{Buns}.  In fact, 
 in \cite{Buns1}, it has been shown that $SR$ decomposition of 
  a general matrix
      could  not be carried out by using
   only the above orthogonal and symplectic transformations. 
 An algorithm, named SRDECO, 
      based on these three transformations  was
    derived in \cite{Buns}.

  From linear
algebra point of view, the $SR$ decomposition via  SRDECO
algorithm does not correspond to the  analog of Householder $QR$
 decomposition, since SRDECO involves transformations which are not  elementary rank-one modification of the identity (transvections), see \cite{Art, Golu}.

In \cite{Sal2} a study, based on  linear algebra concepts and   focusing on the construction of 
  the analog of 
 Householder  transformations in a symplectic linear space, has been accomplished. This has led to the second subset of  transformations. Such analog transformations, which are rank-one modification of the identity are called symplectic Householder transformations. Their main features have been established, especially  the mapping problem has been solved. Then, 
the analog of
   Householder $QR$ decomposition in a symplectic linear space has been derived. 
 The algorithm SRSH
    for computing the
   $SR$ decomposition,  using  these symplectic Householder transformations has been  then presented in
   details.
    Unlike  Householder  $QR$ decomposition, the new algorithm SRSH
     involves free parameters and advantages may be
taken from this fact. It has been demonstrated how these parameters can be
determined in an optimal way providing an optimal
version\cite{Sal3} of the algorithm (SROSH).
  The error analysis and computational aspects of this algorithm have been studied  
\cite{Sal4}. Also,  recently, a mathematical and numerical 
 equivalence between modified symplectic Gram-Schmidt and Householder SR algorithms (typically SRSH or SROSH) have been established  in \cite{Sal5}.
  Computational aspects  and   numerical comparisons
    between SGS and  SROSH have clearly showed the superiority of SROSH over SGS and also that 
  SROSH and SRDECO mostly behave quite similarly, except when SRDECO breaks down. In fact, the latter suffers seriously from the eventuality to encounter a fatal breakdown. The algorithm SROSH works well in these cases, and hence seems to be  adequate  to be used in general, or to be an alternative  to cure the breakdowns in SRDECO.

  In order to build a $SR$-algorithm (which is a $QR$-like algorithm) for computing the eigenvalues and eigenvectors of a matrix \cite{Wat}, a reduction of the matrix to an upper $J$-Hessenberg form is crucial. This is due to the fact that the final algorithm we are looking for should have  $O(n^3)$ as complexity. 

In \cite{Buns}, a reduction of a general matrix to an upper $J$-Hessenberg form is presented, using to this aim, the three symplectic transformations of the above  first subset. The algorithm, called JHESS, is based on an adaptation of  SRDECO.

In this paper, we focus on the reduction of a general matrix, to an upper $J$-Hessenberg form, using only the symplectic Householder transformations (the second subset above).  We show how this reduction can be constructed. The new algorithm, which will be called JHSH algorithm, is based on an adaptation of SRSH algorithm.  A variant of JHSH, named JHOSH is then obtained by taking some  optimal choice of the free parameters. The JHOSH is numerically better than JHSH. However, the accuracy may be lost, since the transformations involved in are not necessarily orthogonal. This leads us to derive another variant, based in replacing when possible, each symplectic non-orthogonal transformation by another one, which is  symplectic and orthogonal.  This gives rise to JHMSH algorithm and its variant JHMSH2.

In this work, we restrict ourselves to the construction of such algorithms. Numerical aspects of the new algorithms and new insights on JHESS algorithm (the choice of the free parameters, near breakdowns, breakdowns, prediction of breakdowns, different strategies of curing near breakdowns, ...) will be studied  separately in a forthcoming paper. Nevertheless,  two illustrating numerical examples   are given, showing in particular  the efficiency of  JHMSH and its variant JHMSH2. More precisely, for these examples, the algorithm JHESS encounter a fatal breakdown, and hence fails to provide any $J$-Hessenberg reduction, while our new algorithms JHMSH, JHMSH2, with a slight modification, perform the $J$-Hessenberg reduction, with a very satisfactory precision for both the errors in the factorization and in the loss of $J$-orthogonality. 

 The remainder of this paper is organized as follows. Section
 2, is devoted to the necessary preliminaries. In the section 3, we  show how we obtain the method of reducing a general matrix to an upper $J$-Hessenberg, based only on the symplectic Householder transformations.  Also, we present two variants, motivated by the numerical stability. 
  Numerical experiments and comparisons between JHESS and the new  JHMSH  are given. We conclude in the section 4. 
\section{Preliminaries}

 Let $J_{2n}$ (or simply $J$) be the $2n$-by-$2n$ real matrix
 \begin{equation}\label{Ji2n}
{ J}_{2n}=\left[\begin{array}{ll} %
{ 0}_n &{ I}_n\\ %
-{ I}_n & { 0}_n
\end{array}
\right],
\end{equation}
 where $0_n$ and $I_n$ stand respectively  for $n$-by-$n$ null and identity
 matrices.
The linear space  $\mathbb{R}^{2n}$  with the indefinite
skew-symmetric inner product
 \begin{equation}
 (x,y)_J = x^T J y
 \end{equation}
  is called symplectic. 
For 
$x,\;y \in\mathbb{R}^{2n},\mbox{ the orthogonality }x \perp' y$ stands for $(x,y)_J=0.$
   The symplectic adjoint $x^J$ of a vector  $x$,  is defined by
  \begin{equation}\label{deff}
  x^J = x^T J.
  \end{equation}
    The symplectic  adjoint of
    $M \in \mathbb{R}^{2n \times 2k}$ is defined  by
 \begin{equation}\label{def3}
 M^{J}=J_{2k}^T M^{T}J_{2n}.
 \end{equation}
 A matrix $S \in \mathbb{R}^{2n \times 2k}$ is called symplectic
 if
 \begin{equation}\label{def2}
 S^{J} S= I_{2k}.
 \end{equation}
 The symplectic  group (multiplicative group of square symplectic matrices) is denoted
 $\mathbb{S}.$
\noindent
A  transformation $T$  given by

\begin{equation}\label{sytr1}
 T=
I+ cvv^J \mbox{ where } c \in \mathbb{R},\; \; v \in
\mathbb{R}^{\nu}\;\; \mbox{ (with } \nu \mbox{ even),}
\end{equation}
 is called symplectic Householder transformation \cite{Sal2}. It satisfies
  \begin{equation}\label{jtra}
  T^J=I - c v
 v^J.
 \end{equation}
  The vector $v$ is called the direction of $T.$

For $x,\;y \in \mathbb{R}^{2n},\;$
  there exists  a symplectic Householder transformation $T$
   such that $Tx=y$
 if $x=y$ or $x^Jy\neq0.$ When $x^Jy\neq0,$
  $T$  is given by
 $$T=I-\frac{1}{x^J y}(y-x)(y-x)^J.$$
Moreover,
 each non null vector $x$ can be mapped onto any non null vector $y$
 by a product of at  most two symplectic Householder transformations \cite{Sal2}.
 Symplectic Householder transformations are rotations, i.e.
$det(T)=1$ and
 the  symplectic group $\mathbb{S}$ is generated by
 symplectic Householder transformations.
We recall that a matrix $H=
 \left[
  \begin{array}{ll}
  H_{11} &H_{12} \\
  H_{21} & H_{22}
  \end{array}
  \right]
   \in  \mathbb{R}^{2n\times 2n},$ is upper $J$-Hessenberg when  $H_{11},\,H_{21},\,H_{22}$ are upper triangular and $H_{12}$ is upper Hessenberg. $H$ is called unreduced when $H_{21}$ is nonsingular and the Hessenberg $H_{12}$ is unreduced, i.e. the entries  of the subdiagonal are all nonzero.
\section{Upper $J$-Hessenberg reduction via symplectic Householder transformations}
\subsection{Toward the algorithm}

  Let $\{e_1,\ldots,e_{2n}\}$ be the canonical basis of $\mathbb{R}^{2n
  }$,  $a\in \mathbb{R}^{2n
  }$ and $\rho,\;\mu,\; \nu$ be arbitrary scalars.
  We seek for symplectic Householder
  transformations $T_1$ and $T_2$ such that
   \begin{equation}\label{algg}
  T_1(a)= \rho e_1,
  \end{equation}
 and 
 \begin{equation}\label{algg1}
  T_2(e_1)=e_1,\;T_2(a)=\mu e_1 + \nu e_{n+1}.
  \end{equation}
   The fact that $T_2 $ is a symplectic isometry yields the necessary
   condition
    \begin{equation}\label{exxx}
     (T_2 (a))^J (T_2 (e_1))=a^J e_1,
    \end{equation}
 which implies $\nu=a{(n+1)}$ and $\mu$ arbitrary. We get
 \begin{thm}\label{algth}
  Let $\rho,\; \mu$ be arbitrary scalars and $\nu =a{(n+1)}.$
  Setting
  $$\displaystyle{c_1=-\frac{1}{ \rho a^J e_1},\; v_1=\rho e_1
  -a,\;c_2=-\frac{1}{a^J (\mu e_1 + \nu e_{n+1})},\, 
  v_2=\mu e_1 + \nu e_{n+1} -a},$$
   then
  \begin{equation}\label{hsel}
  T_1=I+c_1 v_1 v_1^J (\mbox{ respectively }T_2=I+c_2 v_2 v_2^J)\mbox{  satisfy (\ref{algg}) (respectively( \ref{algg1}))}.\;\
  \end{equation}
  
 \end{thm}
\begin{rmk}
Since the $n+1$th component of $v_2$ is zero, $T_2$ keeps the $n+1$th component of $T_2 x$ unchanged, for any $x \in \mathbb{R}^{2n}.$ More on the properties of such transformations $T_1$ or $T_2$ can be found in \cite{Sal3, Sal4}.
\end{rmk}
We also need the following
 \begin{thm}\label{algth2}
Let $v \in \mathbb{R}^{2n}$, with the partition $v=[0^T,u^T,0^T,w^T]^T$, where $[u,w] \in \mathbb{R}^{(n-i)\times 2},$ for a given $\;\;1\leq i \leq n-1$  and set $\tilde{v}=[u^T, w^T]^T$. Consider the symplectic transformations $T= I + c v v^J$ and $ \tilde{T}=I+c \tilde{v} \tilde{v}^J.$ We have \\
$\forall \alpha \in  \mathbb{R}^{i},\;\forall \beta \in \mathbb{R}^{i},\;
\forall x \in  \mathbb{R}^{n-i},\;\forall y \in \mathbb{R}^{n-i},$
$$T[\alpha^T, x^T, \beta^T, y^T]^T=[\alpha^T, x'^T, \beta^T, y'^T]^T, \mbox{ with } [x'^T,y'^T]^T = \tilde{T}[x^T,y^T]^T.$$

 \end{thm}
\begin{pr}
We have $v^J [\alpha^T, x^T, \beta^T, y^T]^T=u^Ty-w^Tx=[u^T w^T]J[x^T y^T]^T.$  Then $T[\alpha^T, x^T, \beta^T, y^T]^T=[\alpha^T, x^T, \beta^T, y^T]^T+ c [0^T,u^T,0^T,w^T]^T [u^T w^T]J[x^T y^T]^T.$
 We check easily 
$\left[\begin{array}{l}
x'\\
y'
\end{array}
\right] =
\left[\begin{array}{l}
x\\
y
\end{array}
\right]
 +c\left[\begin{array}{l}
u\\
w
\end{array}
\right][u^T w^T]J
\left[\begin{array}{l}
x\\
y
\end{array}
\right]=\tilde{T}\left[\begin{array}{l}
x\\
y
\end{array}
\right],$  and \\ $T[\alpha^T, 0^T, \beta^T, 0^T]^T=[\alpha^T, 0^T, \beta^T, 0^T]^T.$  
\end{pr}
Note that the Theorem \ref{algth2} remains valid if one takes $T^J$ instead of $T.$
This result, with Theorem \ref{algth}, constitute the main tool on which the $SR$ factorization (based on symplectic Householder transformations) is constructed.  We will adapt  this tool for  reducing a general matrix to an upper $J$-Hessenberg form,   based on these symplectic Householder transformations.
\subsection{The $J$-Hessenberg reduction : the JHSH algorithm}
   We explain here the
  steps of the algorithm by illustrating  the general pattern. Let $A=[a_1,\ldots,a_n,a_{n+1},\ldots,a_{2n}]\in \mathbb{R}^{2n
  \times 2n}$ be a given matrix and set $A^{(0)}=A.$ We will use the notation $A_{(i_1:i_2,j_1:j_2)}$  to denote the submatrix obtained from the matrix $A$ by deleting all rows and columns except rows $i_1$ until $i_2$ and columns $j_1$ until $j_2.$
\\
{\bf 1.} Choose a symplectic Householder transformation $H_1$  (i.e. $c_1\in\mathbb{R}$ and $v_1 \in \mathbb{R}^{2n}$),  with $H_1 e_1=e_1,$   to zero out entries 2 through $n$ and entries  $n+2$ through $2n$ of the first column of $A$. The vector $e_1$ stands for the first canonical vector of  $\mathbb{R}^{2n}.$
 The transformation $H_1$ corresponds to the transformation $T_2$,  given  in  Theorem \ref{algth}.
 Set $v_1$ 
the direction vector   of  $H_1.$  Since $H_1 e_1=e_1,$ we obtain    $v_1^Je_1=v_1^T J e_1=0$. Thus the $n+1$th component of $v_1$  is zero.   It follows that for any vector $x$,  the $n+1$th component of $H_1 x$  remains unchanged. The direction $v_1$ of $H_1$  
 is given by $v_1= A^{(1)}_{(1,1)}e_{1}+a_1{(n+1)}e_{n+1}-a_1,$  where $ A^{(1)}_{(1,1)}$ is an arbitrary given  scalar.  Notice that   we have also $H_1^Je_1=e_1,$    
and hence the first column of $H_1$ and $H_1^J$ is $e_1.$ 
  Thus, multiplying $A^{(0)}$ on the left by $H_1$ 
   leaves unchanged the $n+1$th row and creates the desired zeros in the first column. We get 
  $$A'^{(1)}=H_1  A^{(0)} =\left[
  \begin{array}{lll}
  A^{(1)}_{(1,1)} & A'^{(1)}_{(1,2:n) }                 & A'^{(1)}_{(1,n+1:2n)}\\
  0      & A'^{(1)}_{(2:n,2:n)}               & A'^{(1)}_{(2:n,n+1:2n)}\\
  A^{(0)}_{(n+1,1)}      & A^{(0)}_{(n+1,2:n)}                & A^{(0)}_{(n+1, n+1:2n)} \\
  0    & A'^{(1)}_{(n+2:2n,2:n)}                & A'^{(1)}_{(n+2:2n, n+1:2n)} 
  \end{array}
  \right].$$
    The step involves the free parameter $ A^{(1)}_{(1,1)}.$
    
     Multiplying $H_1A^{(0)}$ on the right by $H_1^J$ 
     leaves the first column of $H_1A^{(0)} H_1^J$ unchanged, and we obtain
  $$A^{(1)}=H_1  A^{(0)}H_1^J =\left[
  \begin{array}{lll}
  A^{(1)}_{(1,1)} & A^{(1)}_{(1,2:n) }                 & A^{(1)}_{(1,n+1:2n)}\\
  0      & A^{(1)}_{(2:n,2:n)}               & A^{(1)}_{(2:n,n+1:2n)}\\
  A^{(0)}_{(n+1,1)}      & A^{(1)}_{(n+1,2:n)}                & A^{(1)}_{(n+1, n+1:2n)} \\
  0    & A^{(1)}_{(n+2:2n,2:n)}                & A^{(1)}_{(n+2:2n, n+1:2n)} 
  \end{array}
  \right].$$
   The next step consists in choosing a symplectic Householder ${H}_2$ to zero out the entries 3 through n, the entries $n+2$ through $2n$ of the $n+1$th column of $A^{(1)}.$ To do this,  let $\tilde{A}^{(1)}=\left[\begin{array}{ll}
   A^{(1)}_{(2:n,2:n)} &  A^{(1)}_{(2:n,n+1:2n)}\\
   A^{(1)}_{n+2:2n,2:n}&   A^{(1)}_{n+2:2n,n+1:2n}
   \end{array}
   \right]$  be the 
    the matrix obtained from ${A}^{(1)}$ by deleting the first column and the first and the $n+1$th rows. 
And let 
    ${A}^{(2)}_{(2,n+1)} \neq 0$ be an 
arbitrary  given scalar.
   We apply   $\tilde{H}_2=I_{2n-2}+c_2 \tilde{v}_{2}
 {\tilde{v}_{2}}^J$  
  given by   Theorem \ref{algth},
 with $\tilde{v}_{2}=
   \left[
  \begin{array}{l}
  u_2\\
  \hline
  w_2
  \end{array}
  \right]={A}^{(2)}_{(2,n+1)}e_1 -\tilde{A}^{(1)}(:,n)\in  \mathbb{R}^{2n-2},\;u_2\in  \mathbb{R}^{n-1},\;w_2 \in  \mathbb{R}^{n-1}
 ,$  where $e_1$ stands for the first canonical vector of $\mathbb{R}^{2n-2}.$ We obtain
  $$
\tilde{A}'^{(2)}=  \tilde{H}_2 \tilde{A}^{(1)} =\left[
  \begin{array}{lll}
   {A}'^{(2)}_{(2,2:n)} & {A}^{(2)}_{(2,n+1)}& {A}'^{(2)}_{(2,n+2:2n)}\\
   {A}'^{(2)}_{(3:n,2:n)} & 0 & {A}'^{(2)}_{(3:n,n+2:2n)}\\
      {A}'^{(2)}_{(n+2:2n,2:n)} & 0 & {A}'^{(2)}_{(n+2:2n,n+2:2n)}
  \end{array}
  \right].$$
 The transformation $\tilde{H}_2 $ corresponds to the choice $T_1$ in Theorem \ref{algth}.
    Setting  $H_2= I_{2n}+c_2 {v}_{2}
 {{v}_{2}}^J
 ,$ with $
 {v}_{2} 
  =
   \left[
  \begin{array}{l}
  0 \\
  u_2\\
  \hline
  0 \\
  w_2 
  \end{array}
  \right] \in  \mathbb{R}^{2n}$
    then $H_2$ is a symplectic Householder transformation. Using Theorem \ref{algth2},  we get 
 $$
{A}'^{(2)}=  {H}_2 {A}^{(1)} =\left[
  \begin{array}{llll}
  A^{(1)}_{(1,1)} &  {A}^{(1)}_{(1,2:n)} & {A}^{(1)}_{(1,n+1)}& {A}^{(1)}_{(1,n+2:2n)}\\
  0 &  {A}'^{(2)}_{(2,2:n)} & {A}^{(2)}_{(2,n+1)}& {A}'^{(2)}_{(2,n+2:2n)}\\
0  & {A}'^{(2)}_{(3:n,2:n)} & 0 & {A}'^{(2)}_{(3:n,n+2:2n)}\\
  A^{(0)}_{(n+1,1)}  & {A}^{(1)}_{(n+1,2:n)} &  {A}^{(1)}_{(n+1,n+1)}& {A}^{(1)}_{(n+1,n+2:2n)} \\
0  &  {A}'^{(2)}_{(n+2:2n,2:n)}& 0 & { A}'^{(2)}_{(n+2:2n,n+2:2n)}
  \end{array}
  \right].$$
$H_2$ leaves the first and the $n+1$ th rows of ${H}_2 {A}^{(1)}$ unchanged. It leaves the first column of ${H}_2 {A}^{(1)}$ unchanged, and creates the desired zeros in the column $n+1.$

 The multiplication of  ${H}_2 {A}^{(1)}$  on the right by  ${H}^J_2$ leaves the first and the $n+1$th columns of $ {H}_2 {A}^{(1)}{H}^J_2$ unchanged. We obtain

  $$
{A}^{(2)}=  {H}_2{A}^{(1)} H^J_2 =\left[
  \begin{array}{llll}
  A^{(1)}_{(1,1)} &  {A}^{(2)}_{(1,2:n)} & {A}^{(1)}_{(1,n+1)}& {A}^{(2)}_{(1,n+2:2n)}\\
  0 &  {A}^{(2)}_{(2,2:n)} & {A}^{(2)}_{(2,n+1)}& {A}^{(2)}_{(2,n+2:2n)}\\
0  & {A}^{(2)}_{(3:n,2:n)} & 0 & {A}^{(2)}_{(3:n,n+2:2n)}\\
  A^{(0)}_{(n+1,1)}  & {A}^{(2)}_{(n+1,2:n)} &  {A}^{(1)}_{(n+1,n+1)}& {A}^{(2)}_{(n+1,n+2:2n)} \\
0  &  {A}^{(2)}_{(n+2:2n,2:n)}& 0 & { A}^{(2)}_{(n+2:2n,n+2:2n)}
  \end{array}
  \right].$$ It is worth noting that $H_2 e_1=e_1$ and $H_2e_{n+1}=e_{n+1}.$ Thus the first column (respectively  the $n+1$th column) of $H_2$ and $H_2^J$ is $e_1$ (respectively $e_{n+1}$).

  In the next step, we want to zero out the entries  3 through $n$ and $n+3$ through $2n$ of  the second column of $A^{(2)}$ and  the entries  4 through $n$ and $n+3$ through $2n$ of the column $n+2$  of $A^{(2)}.$ Let  $\tilde{A}^{(2)} $ be the matrix obtained from ${A}^{(2)} $ by deleting the first, the $n+1$th rows, and the corresponding columns, ie. 
 $\tilde{A}^{(2)} = \left[
\begin{array}{ll}
{A}^{(2)}_{(2:n,2:n)} & {A}^{(2)}_{(2:n,n+2:2n)}\\
{A}^{(2)}_{(n+2:2n,2:n)} & {A}^{(2)}_{(n+2:2n,n+2:2n)}
\end{array}
\right].$\\
\\
{\bf 2.} We apply now exactly the same two steps of {\bf 1.},  to the new size reduced matrix  $\tilde{A}^{(2)}.$ In other words, we choose a symplectic Householder transformation $\tilde{H}_3$, which means to compute 
a vector $\tilde{v}_3=[u_3^T,w_3^T]^T$ with $u_3\in \mathbb{R}^{n-1},\;\;w_3\in \mathbb{R}^{n-1}$ and a real $c_3$ such that  $\tilde{H}_3 =I+c_3\tilde{ v}_3\tilde{ v}_3^J$  zero out the entries 2 through $n-1$ and the entries $n$ through $2n-2$ of the first column of $\tilde{A}^{(2)}$   with 
 $\tilde{H}_3e_1=e_1 \in \mathbb{R}^{2n-2}$.  The transformation $\tilde{H_3}$ corresponds to the transformation $T_2$,   in  Theorem \ref{algth}.  The direction vector $\tilde{v}_3$ of $\tilde{H_3}$ is given by $\tilde{v}_3=A^{(3)}_{(2,2)}e_{1}+\tilde{A}^{(2)}(n,1)e_{n}-\tilde{A}^{(2)}(:,1),$ where 
 $ A^{(3)}_{(2,2)}$ is an arbitrary non zero scalar. $\tilde{H_3}$ leaves unchanged the $n$th row of $\tilde{H_3}\tilde{A}^{(2)}.$  We get
$$
\tilde{A}'^{(3)}= \tilde{H}_3  \tilde{A}^{(2)}=
\left[
\begin{array}{lll}
{A}^{(3)}_{(2,2)} & {A}'^{(3)}_{(2,3:n)} & {A}'^{(3)}_{(2,n+2: 2n)}\\
0 & {A}'^{(3)}_{(3:n,3:n)} & {A}'^{(3)}_{(3:n,n+2: 2n)}\\
{A}^{(2)}_{(n+2,2)} & {A}^{(2)}_{(n+2,3:n)} & {A}^{(2)}_{(n+2,n+2: 2n)}\\
0& {A}'^{(3)}_{(n+3:2n,3:n)} & {A}'^{(3)}_{(n+3:2n,n+2: 2n)}\\
\end{array}
\right].
$$
Remark that the $n$th component of $\tilde{v}_3$ is zero. Take now $v_3=[0\; u_3^T|0\; w_3^T]^T$ 
%
and set $H_3= I + c_3 v_3 v_3^J$. Then  $H_3$ is obviously a symplectic Householder transformation of order $2n$. The components 1, $n+1$ and $n+2$ of $v_3$ are equal to zero.  Thus $H_3$ leaves the rows 1,  $n+1$ and $n+2$ of $H_3 A^{(2)}$ unchanged and satisfy 
$H_3(e_1)=e_1$, $H_3 e_2=e_2$ and $H_3 e_{n+1} = e_{n+1}.$  Thus 

$H_3$ leaves the first and the $n+1$th columns of $H_3 A^{(2)}$ unchanged and zero out the entries 3 through $n$ and the entries $n+3$ through $2n$ of the second column.

 We have 
$$
{A}'^{(3)}= {H}_3 {A}^{(2)}=
\left[
  \begin{array}{lllll}
  A^{(1)}_{(1,1)} &     {A}^{(2)}_{(1,2)}    &{A}^{(2)}_{(1,3:n)} & {A}^{(1)}_{(1,n+1)}& {A}^{(2)}_{(1,n+2:2n)}\\
  0 &                               {A}^{(3)}_{(2,2)}   &{A}'^{(3)}_{(2,3:n)} & {A}^{(2)}_{(2,n+1)}& {A'}^{(3)}_{(2,n+2:2n)}\\
0  &                              0      &{A}'^{(3)}_{(3:n,3:n)} & 0 & {A'}^{(3)}_{(3:n,n+2:2n)}\\
  A^{(0)}_{(n+1,1)}  &         {A}^{(2)}_{(n+1,2)}   & {A}^{(2)}_{(n+1,3:n)} &  {A}^{(1)}_{(n+1,n+1)}& {A}^{(2)}_{(n+1,n+2:2n)} \\
0  &                                  {A}^{(2)}_{(n+2,2)}           & {A}^{(2)}_{(n+2,3:n)}& 0 &{ A}^{(2)}_{(n+2,n+2:2n)}\\
0  &                      0                      & {A}'^{(3)}_{(n+3:2n,3:n)}& 0 & { A'}^{(3)}_{(n+3:2n,n+2:2n)}

  \end{array}
  \right].
$$

The transformation $H^J_3$ leaves the column 1, 2 and $n+1$  of 
 ${H}_3 {A}^{(2)}H_3^J$ unchanged 
since $H^J_3(e_1)=e_1,\;H^J_3 e_2=e_2$ and $H^J_3 e_{n+1} = e_{n+1}.$ 
We get 
$$
{A}^{(3)}= {H}_3  {A}^{(2)}H_3^J=
\left[
  \begin{array}{lllll}
  A^{(1)}_{(1,1)} &     {A}^{(2)}_{(1,2)}    &{A}^{(3)}_{(1,3:n)} & {A}^{(1)}_{(1,n+1)}& {A}^{(3)}_{(1,n+2:2n)}\\
  0 &                               {A}^{(3)}_{(2,2)}   &{A}^{(3)}_{(2,3:n)} & {A}^{(2)}_{(2,n+1)}& {A}^{(3)}_{(2,n+2:2n)}\\
0  &                              0      &{A}^{(3)}_{(3:n,3:n)} & 0 & {A}^{(3)}_{(3:n,n+2:2n)}\\
  A^{(0)}_{(n+1,1)}  &         {A}^{(2)}_{(n+1,2)}   & {A}^{(3)}_{(n+1,3:n)} &  {A}^{(1)}_{(n+1,n+1)}& {A}^{(3)}_{(n+1,n+2:2n)} \\
0  &                                  {A}^{(2)}_{(n+2,2)}           & {A}^{(3)}_{(n+2,3:n)}& 0 & { A}^{(3)}_{(n+2,n+2:2n)}\\
0  &                      0                      & {A}^{(3)}_{(n+3:2n,3:n)}& 0 & { A}^{(3)}_{(n+3:2n,n+2:2n)}

  \end{array}
  \right].
$$
\\
Now, deleting the rows $1,\, 2,\, n+1,\,n+2$ and the columns $1,\, 2,\,n+1$ of $A^{(3) }$ and setting $\tilde{A}^{(3)}=\left[
\begin{array}{ll}
A^{(3)}_{(3:n,3:n)} &A^{(3)}_{(3:n,n+2:2n)}  \\
A^{(3)}_{(n+3:2n,3:n)} &A^{(3)}_{(n+3:2n,n+2:2n)}  
\end{array}
\right]
,$ we find  $c_4 \in \mathbb{R}$ and 
$\tilde{v}_4=\left[
\begin{array}{l}
u_4\\
\hline
w_4
\end{array}
\right],$
  with $u_4 \in \mathbb{R}^{n-2}$ and $w_4 \in \mathbb{R}^{n-2} $ such that the action of 
$\tilde{H}_4=I + c_4 \tilde{v}_4\tilde{v_4}^J $ gives 
$$\tilde{A}'^{(4)}=\tilde{H}_4\tilde{A}^{(3)}=\left[
\begin{array}{lll}
{A}'^{(4)}_{(3,3:n)} & {A}^{(4)}_{(3,n+2)}  &{A}'^{(4)}_{(3,n+3:2n)}  \\
{A}'^{(4)}_{(4:n,3:n)} & 0 &{A}'^{(4)}_{(4:n,n+3:2n)}  \\
{A}'^{(4)}_{(n+3:2n,3:n)}& 0&{A}'^{(4)}_{(n+3:2n,n+3:2n)}  
\end{array}
\right].
$$
 The coefficient  ${A}^{(4)}_{(3,n+2)} $ is an arbitrary chosen scalar.
Taking  $v_4=[0\;0\;u_4^T|0\;0\;w_4^T]^T$
%
then the  transformation  $H_4=I+ c_4 v_4 v^J_4$  leaves unchanged the rows 1, 2, $n+1,\,n+2$  and columns 1, 2, and $n+1$  of $A'^{(4)}=H_4 A^{(3)} $ and creates the desired zeros in the column $n+2.$ We obtain 
$$A'^{(4)}=
\left[
  \begin{array}{llllll}
  A^{(1)}_{(1,1)} &     {A}^{(2)}_{(1,2)}    &{A}^{(3)}_{(1,3:n)} & {A}^{(1)}_{(1,n+1)}&{A}^{(3)}_{(1,n+2)} &{A}^{(3)}_{(1,n+3:2n)}\\
  0 &                               {A}^{(3)}_{(2,2)}   &{A}^{(3)}_{(2,3:n)} & {A}^{(2)}_{(2,n+1)}&{A}^{(3)}_{(2,n+2)} &{A}^{(3)}_{(2,n+3:2n)}\\

0  &                              0      &{A'}^{(4)}_{(3,3:n)} & 0 &  {A}^{(4)}_{(3,n+2)} &{A'}^{(4)}_{(3,n+3:2n)}\\
0  &                              0      &{A'}^{(4)}_{(4:n,3:n)} & 0 &  0 &{A'}^{(4)}_{(4:n,n+3:2n)}\\
  A^{(0)}_{(n+1,1)}  &         {A}^{(2)}_{(n+1,2)}   & {A}^{(3)}_{(n+1,3:n)} &  {A}^{(1)}_{(n+1,n+1)}&       
{A}^{(3)}_{(n+1,n+2)} &{A}^{(3)}_{(n+1,n+3:2n)} \\
0  &                                  {A}^{(2)}_{(n+2,2)}           & {A}^{(3)}_{(n+2,3:n)}& 0 &   { A}^{(3)}_{(n+2,n+2)}  &{ A}^{(3)}_{(n+2,n+3:2n)}\\
0  &                      0                      & {A'}^{(4)}_{(n+3:2n,3:n)}& 0 &  0  &{ A'}^{(4)}_{(n+3:2n,n+3:2n)}
\end{array}
  \right].
$$
$H^J_4$ leaves unchanged the first, the second, the $n+1$, $n+2$ columns of $A^{(4)}=H_4 A^{(3)} H^J_4$ since $H^J_4 (e_i)=e_i$ for $i=1,\,2\,,n+1,\,n+2.$ Hence, we  get
$$A^{(4)}=
\left[
\begin{array}{llllll}
A^{(1)}_{(1,1)} &     {A}^{(2)}_{(1,2)}    &{A}^{(4)}_{(1,3:n)} & {A}^{(1)}_{(1,n+1)}&{A}^{(3)}_{(1,n+2)} &{A}^{(4)}_{(1,n+3:2n)}\\
0 &                               {A}^{(3)}_{(2,2)}   &{A}^{(4)}_{(2,3:n)} & {A}^{(2)}_{(2,n+1)}&{A}^{(3)}_{(2,n+2)} &{A}^{(4)}_{(2,n+3:2n)}\\

0  &                              0      &{A}^{(4)}_{(3,3:n)} & 0 &  {A}^{(4)}_{(3,n+2)} &{A}^{(4)}_{(3,n+3:2n)}\\
0  &                              0      &{A}^{(4)}_{(4:n,3:n)} & 0 &  0 &{A}^{(4)}_{(4:n,n+3:2n)}\\
A^{(0)}_{(n+1,1)}  &         {A}^{(2)}_{(n+1,2)}   & {A}^{(4)}_{(n+1,3:n)} &  {A}^{(1)}_{(n+1,n+1)}&       
{A}^{(3)}_{(n+1,n+2)} &{A}^{(4)}_{(n+1,n+3:2n)} \\
0  &                                  {A}^{(2)}_{(n+2,2)}           & {A}^{(4)}_{(n+2,3:n)}& 0 &   { A}^{(3)}_{(n+2,n+2)}  &{ A}^{(4)}_{(n+2,n+3:2n)}\\
0  &                      0                      & {A}^{(4)}_{(n+3:2n,3:n)}& 0 &  0  &{ A}^{(4)}_{(n+3:2n,n+3:2n)}
\end{array}
\right].$$
{\bf 3.}   The $j$th step is now clear. It  involves  two sub-steps. The first  consists  in finding  $H_{2j-1}$ , i.e. the scalar $c_{2j-1}$ and the vector $v_{2j-1}$ such that $H_{2j-1}=I + c_{2j-1}v_{2j-1}v_{2j-1}^J$ leaves the rows $1,\ldots, j-1,$ the rows  $n+1,\ldots,n+j,$  
the columns $1,\ldots,j-1,$ and the columns  $n+1,\ldots,n+j-1$ of $H_{2j-1}A^{(2j-2)}$  unchanged and zero out the entries $j+1$ through $n$ and the entries $n+j+1$ through $2n$ of the $j$th column. The vector $v_{2j-1}\in \mathbb{R}^{2n}$ has the structure 
$v_{2j-1}=[0^T,u_{2j-1}^T,0^T,w_{2j-1}^T]^T,$ with $u_{2j-1}\in \mathbb{R}^{n-j+1},\;w_{2j-1}\in \mathbb{R}^{n-j+1}.$ The first component of $w_{2j-1}$ is zero.Thus $H_{2j-1}e_i=e_i$ for $i=1,\ldots,j$ and for $i=n+1,\ldots, n+j-1.$ The $j$th column  $H_{2j-1} A^{(2j-2)}(:,j)$ is transformed as follows 
$$H_{2j-1}A^{(2j-2)}(:,j)=\left[
\begin{array}{l}
A^{(2j-2)}(1:j-1,j)\\
A^{(2j-1)}(j,j)\\
0\\
A^{(2j-2)}(n+1:n+j,j)\\
0
\end{array}
\right]
\begin{array}{l}
\left\{j-1\right\} \\
\left\{1\right\}\\
\{n-j\}\\
\left\{j\right\}\\
\left\{n-j\right\}
\end{array}.
$$
The entry $A^{(2j-1)}(j,j)$ is a free parameter. 

The multiplication of $H_{2j-1}A^{(2j-2)}$ on the right by $H_{2j-1}^J$
leaves the columns 
$1,\ldots,j,$ and the columns $n+1,\ldots,n+j-1,$ of $H_{2j-1}A^{(2j-2)}H_{2j-1}^J$  unchanged. 
 The coefficient $c_{2j-1}$, the vector $v_{2j-1}$ and hence  the symplectic transformation $H_{2j-1}$  are simply and explicitly  given by Theorem \ref{algth}. 
 The matrix $A^{(2j-1)}=H_{2j-1}A^{(2j-2)}H^J_{2j-1}$ has the desired form.
  Let us set  $\tilde{H}_{2j-1}=I + c_{2j-1}\tilde{v}_{2j-1}\tilde{v}_{2j-1}^J$,  
  $\tilde{v}_{2j-1}=[u_{2j-1}^T,w_{2j-1}^T]^T,$ where 
   $  [u_{2j-1},w_{2j-1}]\in \mathbb{R}^{\alpha_j\times 2}$ , with $\alpha_j= n-j+1$ 
    and $\tilde{A}^{(2j-2)}(:,j)$ the $j$th column of $\tilde{A}^{(2j-2)}$ obtained from ${A}^{(2j-2)}(:,j)$ by deleting the rows $1,\ldots,j-1$ and rows $n+1,\ldots,n+j-1.$ We obviously obtain $\tilde{H}_{2j-1} \tilde{A}^{(2j-2)}(:,j)= A^{(2j-1)}(j,j) e_1 + A^{(2j-2)}(n+j,j)e_{\alpha_j + 1}.$
     Here $e_1$ and $e_{\alpha_j+1}$ denote the first and the $\alpha_j +1$th canonical vectors of $\mathbb{R}^{2\alpha_j}.$

 In a similar way, the second sub-step consists  in finding  $H_{2j}$ , i.e. the scalar $c_{2j}$ and the vector $v_{2j}$ such that $H_{2j}=I + c_{2j}v_{2j}v_{2j}^J$ leaves the rows $1,\ldots, j,$ the rows  $n+1,\ldots,n+j,$  
 the columns $1,\ldots,j,$ and the columns  $n+1,\ldots,n+j-1$ of $H_{2j}A^{(2j-1)}$  unchanged and zero out the entries $j+2$ through $n$ and the entries $n+j+1$ through $2n$ of the $n+j$th column. The vector $v_{2j}\in \mathbb{R}^{2n}$ has the structure 
 $v_{2j}=[0^T,u_{2j}^T,0^T,w_{2j}^T]^T,$ with $u_{2j}\in \mathbb{R}^{n-j},\;w_{2j}\in \mathbb{R}^{n-j}.$ 
 Thus $H_{2j}e_i=e_i$ for $i=1,\ldots,j$ and for $i=n+1,\ldots, n+j.$ The $n+j$th column of $H_{2j} A^{(2j-1)}(:,n+j)$ is transformed as follows 
 $$H_{2j}A^{(2j-1)}(:,n+j)=\left[
 \begin{array}{l}
 A^{(2j-1)}(1:j,n+j)\\
 A^{(2j)}(j+1,n+j)\\
 0\\
 A^{(2j-1)}(n+1:n+j,n+j)\\
 0
 \end{array}
 \right]
 \begin{array}{l}
 \left\{j\right\} \\
 \left\{1\right\}\\
 \{n-j-1\}\\
 \left\{j\right\}\\
 \left\{n-j\right\}
 \end{array}.
 $$
 The entry $A^{(2j)}(j+1,n+j)$ is a free parameter. 
 
 The multiplication of $H_{2j}A^{(2j-1))}$ on the right by $H_{2j}^J$
 leaves the columns 
 $1,\ldots,j,$ and the columns $n+1,\ldots,n+j,$ of $H_{2j}A^{(2j-1)}H_{2j}^J$  unchanged. 
  The coefficient $c_{2j}$, the vector $v_{2j}$ and hence  the symplectic transformation $H_{2j}$  are  explicitly  given by Theorem \ref{algth}. 
   The matrix $A^{(2j)}=H_{2j}A^{(2j-1)}H^J_{2j}$ has the desired form.
  
  Let us set  $\tilde{H}_{2j}=I + c_{2j}\tilde{v}_{2j}\tilde{v}_{2j}^J$,  with 
  $\tilde{v}_{2j}=[u_{2j}^T,w_{2j}^T]^T,$ where 
  $  [u_{2j},w_{2j}]\in \mathbb{R}^{\beta_j\times 2}$ , $\beta_j= n-j$ 
  and $\tilde{A}^{(2j-1)}(:,n+j)$ the $n+j$th column of $\tilde{A}^{(2j-1)}$ obtained from ${A}^{(2j-1)}(:,n+j)$ by deleting the rows $1,\ldots,j$ and rows $n+1,\ldots,n+j.$ We obviously obtain $\tilde{H}_{2j} \tilde{A}^{(2j-1)}(:,n+j)= A^{(2j)}(j+1,n+j) e_1.$
  Here $e_1$ denotes the first canonical vector of $\mathbb{R}^{2\beta_j.}$

  Thus,  it is worth noting that each step $j$ involves two free parameters 
  $A^{(2j-1)}(j,j)$  and $A^{(2j)}(j+1,n+j),$ and that these parameters are located as highlighted above, in the corresponding symplectic Householder transformations $H_{2j-1}$ and $H_{2j}$  (or equivalently 
   $\tilde{H}_{2j-1}$ and $\tilde{H}_{2j}$ ).

  At the last step (the $n-1$th step), we obtain\\
  $ H_{2n-2} \ldots H_2 H_1 A ( H_{2n-2} \ldots H_2 H_1 )^J=  \left[
  \begin{array}{ll}
  H_{11} &H_{12} \\
  H_{21} & H_{22}
  \end{array}
  \right]
  =H \in  \mathbb{R}^{2n\times 2n},$ with $H_{11},\,H_{21},\,H_{22}$
  upper triangular and $H_{12}$  upper Hessenberg. 
    We get $A=S^JHS$ with 
$S=H_{2n-2} \ldots   H_{1}.$
The entries of the diagonal of $H_{11}$ are the free parameters $A^{(2j-1)}(j,j),$ ie. $H_{11}(j,j)=A^{(2j-1)}(j,j)\mbox{ for }j=1,\ldots,n.$  Also, The entries of the sub-diagonal of $H_{12}$ are the free parameters $A^{(2j)}(j+1,n+j),$ ie. $H_{12}(j+1,j)=A^{(2j)}(j+1,n+j)  \mbox{ for }j=1,\ldots,n-1.$
We propose here the  algorithm in its general version, written in pseudo Matlab code,  for computing the reduction of a matrix to the upper $J$-Hessenberg form, via symplectic Householder transformations  (JHSH algorithm).
\begin{alg}
	function [S,H]=JHSH(A)\\
	\indent $twon=size(A(:,1));\;n=twon/2;\;
	S=eye(twon);$\\
	\indent $\mbox{for }  j=1:n-1$\\
	\indent   $J=[zeros(n-j+1),eye(n-j+1);-eye(n-j+1),zeros(n-j+1)];$\\
	\indent    $ro=[j:n,n+j:2n];\;co=[j:n,n+j:2n];$\\
	\indent   $[c,v]=sh2(A(ro,j));$\\
	\indent \% Updating $A:$\\
	\indent    $A(ro,co)=A(ro,co)+c * v * (v' * J * A(ro,co));$\\
	\indent    $A(:,co)=A(:,co)-(A(:,co) * (c * v)) * v' * J;$\\
	\indent\% Updating $S$ (if needed):\\
	\indent    $S(ro,2:end) = S(ro,2:end)+ c * (v * v') * J * S(ro,2:end);$ 
	\indent    $J=[zeros(n-j),eye(n-j);-eye(n-j),zeros(n-j)];$\\
	\indent    $ro=[j+1:n,n+j+1:2n];$\\
	\indent    $[c,v]=sh1(A(ro,n+j));$\\
	\indent \%Updating $A$:\\
	\indent    $A(ro,co)=A(ro,co)+c * v * (v' * J * A(ro,co));$\\    
	\indent   $A(:,co)=A(:,co)-(A(:,co) * (c * v)) * v' * J;$\\
	\indent \%Updating $S$ (if needed):\\
	\indent    $S(ro,2:end) = S(ro,2:end) + c * (v * v') * J * S(ro,2:end);$  
	 \\
	\indent end\\
	\indent end\\
\end{alg}
\begin{alg}\label{sh1}
	function [c, v] = sh1(a)\\
	\indent \%compute $c$ and $v$ such that $T_1 a = \rho e_1,$ \\
	\indent \%$\rho$ is a free parameter, and  $T_1=(eye(twon)+c * v * v' * J);$\\
	\indent $twon=length(a);\;n=twon/2;$\\
	\indent $J=[zeros(n),eye(n);-eye(n),zeros(n)];$\\
	\indent $\mbox{choose }\rho;\;aux = a(1)-\rho;$\\
	\indent if $aux==0$\\
	\indent    $c=0;\;v=zeros(twon,1);\;\%T=eye(twon);$\\    
	\indent elseif $a{(n+1)}==0$\\
	\indent       display('division by zero');\\
	\indent       return\\
	\indent else\\
	\indent    $\displaystyle{v=\frac{a}{aux}};\;\displaystyle{c=\frac{aux^{2}}{\rho \times a{(n+1)}}};\;v(1)=1;$\\
	\indent end\\
	\indent end\\
\end{alg}
\begin{alg}\label{sh2}
	function [c, v] = sh2(a)\\
		\indent \%compute $c$ and $v$ such that $T_2 e_1 = e_1,$ and $T_2 a = \mu e_1 + \nu e_{n+1},$ \\
		\indent \%$\mu$ is a free parameter, and  $T_2=(eye(twon)+c * v * v' * J);$\\
	\indent $twon=length(a);\;n=twon/2;$\\
	\indent $J=[zeros(n),eye(n);-eye(n),zeros(n)];$\\
	\indent if $n==1$\\
	\indent\indent    $v=zeros(twon,1);\;c=0;$   $\%T=eye(twon);$\\
	\indent else\\
\indent choose $\mu;$\\
	\indent       $ \nu =a(n+1);$\\
	\indent        if $ \nu ==0$\\
	\indent           display('division by zero')\\
	\indent           return\\
	\indent       else\\
	\indent           $v=\mu e_1 + \nu e_{n+1} - a,$ 
	 $c=\frac{1}{a(n+1)(a(1)-\mu)};$\\
%
	\indent  end\\
	\indent end
\end{alg}
\subsection{JHOSH, JHMSH algorithms}
 From   an algebraic point of view, JHSH is the analog in the symplectic  case, of the algorithm performing the Hessenberg reduction of a matrix via Householder transformations in the Euclidean case.  Recall that JHSH involves two free parameters at each steps, and the involved symplectic Householder transformations are not orthogonal.  In the sequel, we show how one can take benefit from these free parameters in some optimal way. In order to get an algorithm numerically stable as possible, the free parameters will be chosen so that the symplectic Householder transformations used in the reduction have minimal  norm-2 condition number. The choice of such parameters is as follows \cite{Sal3} :
\begin{thm}\label{algth3}
 Let $\{e_1,\ldots,e_{2n}\}$ be the canonical basis of $\mathbb{R}^{2n
 }$ and  $a\in \mathbb{R}^{2n
}$ given. Take   $\rho=sign(a(1))\norm{a}_2$ and $\mu=a(1)\pm\xi,\; \nu=a(n+1)$  with $\xi=\sqrt{\sum_{i=2, i\neq n+1}^{2n}a(i)^2}.$ 
	Setting
	$$\displaystyle{c_1=-\frac{1}{ \rho a^J e_1},\; v_1=\rho e_1
		-a,\;c_2=-\frac{1}{a^J (\mu e_1 + \nu e_{n+1})},\, 
		v_2=\mu e_1 + \nu e_{n+1} -a},$$
	then
	\begin{equation}\label{hsel}
	T_1=I+c_1 v_1 v_1^J (\mbox{ respectively }T_2=I+c_2 v_2 v_2^J)\mbox{  satisfy (\ref{algg}) (respectively( \ref{algg1}))},\;\
	\end{equation}
with $T_1$ (respectively $T_2$)	 has the minimal norm-2 condition number.
\end{thm}
\begin{pr}
	See \cite{Sal3}.
	\end{pr}
	 For these choices of the free parameters, we refer to $T_1$ (respectively $T_2$) as the first optimal symplectic Householder (osh1) transformation (respectively the second optimal symplectic Householder osh2) transformation.
 This optimal version of JHSH is referred to as JHOSH algorithm and is given as follows :
\begin{alg}
function [S,H]=JHOSH(A)\\
\indent replace in the body of JHSH the sh1 by osh1 and sh2 by osh2. \\
end.
\end{alg}
The pseudo code Matlab of $osh1$ and $osh2$ is a follows
\begin{alg}\label{opts1}
function [c, v] = osh1(a)\\
\indent $twon=length(a);\;n=twon/2;$\\
\indent $J=[zeros(n),eye(n);-eye(n),zeros(n)];$\\
\indent $\rho=sign(a(1))* \Vert a \Vert _{2};\;aux = a(1)-\rho;$\\
\indent if $aux==0$\\
\indent   $c=0;\;v=zeros(twon,1);\;\%T=eye(twon);$\\    
\indent elseif $a{(n+1)}==0$\\
\indent       display('division by zero');\\
\indent        return\\
\indent else\\
\indent    $\displaystyle{v=\frac{a}{aux}};\;\displaystyle{c=\frac{aux^{2}}{\rho * a{(n+1)}}};\;v(1)=1;$\\
\indent   $\%T=(eye(twon)+c * v * v' * J);$\\
\indent end\\
\indent end\\
\end{alg}
\begin{alg}\label{opts2}
function [c, v] = osh2(a)\\
\indent $twon=length(u);\;n=twon/2;$\\
\indent $J=[zeros(n),eye(n);-eye(n),zeros(n)];$\\
\indent if $n==1$\\
\indent    $v=zeros(twon,1);\;c=0;$   $\%T=eye(twon);$\\
\indent else\\
\indent    $I=[2:n,n+2:twon];\;\xi=norm(a(I));$\\
\indent    if $\xi==0$\\
\indent        $v=zeros(twon,1);\;c=0;\;$      $\%T=eye(twon);$\\
\indent    else\\
\indent        $ \nu =a(n+1);$\\
\indent       if $ \nu ==0$\\
\indent            display('division by zero')\\
\indent           return\\
\indent       else\\
\indent          $v=-a/\xi;\;           v (1)=1;\;v(n+1)=0;\;c=\xi/\nu;$\\
\indent            $\%T=(eye(twon)+c * v * v' * J);$\\
\indent        end\\
\indent    end\\
\indent end\\
\indent end
\end{alg}
We have seen that the symplectic Householder transformations used in JHOSH algorithm have minimal norm-2 condition number, and thus numerically,  JHOSH  presents  a significant  advantage over JHSH. However, all these symplectic Householder transformations are not orthogonal. It is well known that it is not possible to handle a $SR$ decomposition  using only transformations which are both symplectic and orthogonal (see \cite{Buns1}). Nevertheless, we will show that half of them (all the transformations $H_{2j}$ above)  may be replaced by specified transformations which are both orthogonal and symplectic. Furthermore, we will show that the two type of orthogonal and symplectic transformations, introduced by Paige et al. \cite{Paig, Vloa}
  can be used to replace the symplectic transformations $H_{2j}$, 
 to zero desired components of a vector.  The first type is
 \begin{equation}\label{reff1}
 H(k,w)=\left(
 \begin{array}{ll}
 \mbox{diag}\displaystyle{(I_{k-1},P)} & 0\\
 0 & \mbox{diag}\displaystyle{ (I_{k-1},P)}
 \end{array}
 \right),
 \end{equation}
 where
 $$P=I-2w w^T/w^Tw, \;\;w \in \mathbb{R}^{n-k+1}.$$
  The transformation 
 $H(k,w)$ is just a direct sum of
 two "ordinary" $n-$by$-n$ Householder matrices \cite{Wilk}. We refer to $H(k,w)$ as Van Loan's Householder transformations. The
 second type is
 \begin{equation}\label{reff2}
 J(k,c,s)=\left(
 \begin{array}{ll}
 C & S\\
 -S & C
 \end{array}
 \right),
 \end{equation}
 where $c^2+s^2=1,$ and 
 $$C=diag(I_{k-1},c,I_{n-k}),$$
 $$S=diag(0_{k-1},s,0_{n-k}).$$
 $J(k,c,s)$ is a  Givens  transformation, which is an "ordinary" 2$n$-by-$2n$ Givens rotation that
 rotates in planes $(k,k+n)$ \cite{Wilk}. We refer to $J(k,c,s)$ as Van Loan's Givens rotation.
 Van Loan's Householder and Givens  transformations  are both
 orthogonal and symplectic. It is worth noting that for $i\neq k $ and $i\neq n+k,$ we have $J(k,c,s)e_i=e_i.$  Also, we have $J(k,c,s)e_k=c e_k -s e_{n+k}$ and $J(k,c,s)e_{n+k}=s e_k +c e_{n+k}.$ Thus, $J(k,c,s)$ leaves unchanged all the  rows of $J(k,c,s)a$ except rows $k$ and $n+k.$ It is obvious also that $H(k,w)e_i=e_i$ for $i=1,\ldots,k-1$ and $i=n+1,\ldots,n+k-1.$  The modification of the even sub-steps of JHOSH (or JHSH) algorithm is as follows.  Let $A=[a_1,\ldots,a_n,a_{n+1},\ldots,a_{2n}]\in \mathbb{R}^{2n
 	\times 2n}$ be a given matrix and set $A^{(0)}=A.$ The first sub-step is obtained by creating the desired zeros in the first column, via the $H_1$ as above. The updated matrix is $A^{(1)}.$ Now, for creating the desired zeros in the column $n+1$ and keeping the first column unchanged, we shall use the Van Loan's transformations, instead of $H_2.$ For $k=n,\ldots,2$, we compute $J(k,c,s)$ such that a zero is created in position $n+k$ in the $n+1$th column of $J(k,c,s) A^{(1)}.$
  The first column as well as the already created zeros in the current $n+1$ column of $A^{(1)}$  remain unchanged. The first and the $n+1$th columns of $J(k,c,s) A^{(1)}$ leave unchanged when the latter is multiplied on the right by  $J(k,c,s)^T.$ The matrix $A^{(1)}$ is then  updated with $A^{(2)}=J(k,c,s)A^{(1)}J(k,c,s)^T.$ So the entries  at  positions $n+2,\ldots,2n$ in  the $n+1$ column of $A^{(2)}$ are zeros. 
  Now, we compute $w$ so that the action of Van Loan's Householder in the product  $H(2,w)A^{(2)}$ creates zeros in the positions $3,\ldots,n$ in the $n+1$ column. The first column of $H(2,w)A^{(2)}$ as well as the already created zeros remain unchanged.   The transformation $H(2,w)$ leaves unchanged the first and the $n+1$ columns of the updated matrix $A^{(2)}=H(2,w)A^{(2)}H(2,w)^T.$ \\
  At the $j$th step, the first sub-step is obtained by creating the desired zeros in the $j$th column , via the $H_{2j-1}$ as in JHOSH. The updated matrix is $A^{(2j-1)}.$ Now, the desired zeros in the column $n+j$ are created by using the Van Loan's givens rotations, instead of $H_{2j}.$ 
   For $k=n,\ldots,j+1$, we compute $J(k,c,s)$ such that a zero is created in position $n+k$ in the $n+j$th column of $J(k,c,s) A^{(2j-1)}.$
  The columns $1,\ldots,j$ and $n+1,\ldots,n+j-1$ as well as the already created zeros in the current $n+j$ column of $A^{(2j-1)}$  remain unchanged. 
  The columns $1,\ldots,j$ and $n+1,\ldots,n+j$
   of $J(k,c,s) A^{(2j-1)}$ leave unchanged when the latter is multiplied on the right by  $J(k,c,s)^T.$ The matrix $A^{(2j-1)}$ is then  updated with $A^{(2j)}=J(k,c,s)A^{(2j-1)}J(k,c,s)^T.$ So the entries  at  positions $n+j+1,\ldots,2n$ in  the $n+j$ column of $A^{(2j)}$ are zeros. 
  Now, we compute $w$ so that the action of Van Loan's Householder in the product  $H(j,w)A^{(2j)}$ creates zeros in the positions $j+2,\ldots,n$ in the $n+j$th column. 
   The columns $1,\ldots,j$ and $n+1,\ldots,n+j-1$ as well as the already created zeros in the current $n+j$ column of $A^{(2j)}$  remain unchanged. 
  $H(j,w)$ leaves unchanged the 
  columns $1,\ldots,j$ and $n+1,\ldots,n+j$
   of the updated matrix $A^{(2j)}=H(j,w)A^{(2j)}H(j,w)^T.$ 
    We obtain the following algorithm
   \begin{alg}
   	function [S,H]=JHMSH(A)\\
   	\indent $twon=size(A(:,1));\;n=twon/2;\;
   	S=eye(twon);$\\
   	\indent $for j=1:n-1$\\
   	\indent\indent    $J=[zeros(n-j+1),eye(n-j+1);-eye(n-j+1),zeros(n-j+1)];$\\
   	\indent\indent    $ro=[j:n,n+j:2n];\;co=[j:n,n+j:2n];$\\
   	\indent\indent    $[c,v]=osh2(A(ro,j));$\\
   	\indent \% Updating $A:$\\
   	\indent\indent    $A(ro,co)=A(ro,co)+c * v * (v' * J * A(ro,co));$\\
   	\indent\indent    $A(:,co)=A(:,co)-(A(:,co) * (c * v)) * v' * J;$\\
   	\indent\% Updating $S$ (if needed):\\
   	\indent\indent    $S(:,co) = S(:,co)- c * (v * v') * J * S(:,co);$ 
   \\
   	\indent \indent for $k=2n:n+j+1,$\\
   	\indent \indent  $[c,s]=vlg(k,A(:,n+j)),\;$\\
   	\indent	\indent \%Updating $A$:\\
     	\indent	\indent $\left[
    	\begin{array}{l}
    	A(k,co) \\
    	A(n+k,co)	
    		\end{array}
    	\right]=
    		\left[
    		\begin{array}{ll}
    			c & s\\
    			-s & c	
    		\end{array}
    		\right]
    		\left[
    		\begin{array}{l}
    			A(k,co) \\
    			 A(n+k,co)	
    		\end{array}
    		\right];$\\
    		\indent	\indent $\left[
    		\begin{array}{ll}
    		A(:,k) &	A(:,n+k)	
    		\end{array}
    		\right]=
  \left[
  \begin{array}{ll}
  A(:,k) &	A(:,n+k)	
  \end{array}
  \right]  		
    	\left[
    		\begin{array}{ll}
    		c & -s\\
    		s & c	
    		\end{array}
    		\right];
    	$\\
   	\indent \%Updating $S$ (if needed):\\
     	\indent	\indent $\left[
     	\begin{array}{ll}
     	S(:,k) &	S(:,n+k)	
     	\end{array}
     	\right]=
     	\left[
     	\begin{array}{ll}
     	S(:,k) &	S(:,n+k)	
     	\end{array}
     	\right]  		
     	\left[
     	\begin{array}{ll}
     	c & -s\\
     	s & c	
     	\end{array}
     	\right];
     	$\\ 	
   	 	\indent\indent end\\
   	\indent \indent if $j\leq n-2$\\
    \indent\indent	[$\beta$,w]=vlh(j+1,A(:,n+j));\\
   	 	\indent	\indent \%Updating $A$:\\
   	 	\indent\indent    $A(j+1:n,co)=A(j+1:n,co)-\beta *w*w'*A(j+1:n,co) $\\
   	 	\indent\indent    $A(j+1+n:2n,co)=A(j+1+n:2n,co)-\beta *w*w'*A(j+1+n:2n,co);$\\	
   	 	\indent\indent    $A(:,j+1:n)=A(:,j+1:n)-\beta *A(:,j+1:n)w*w';$\\	
   	 	\indent\indent    $A(:,n+j+1:2n)=A(:,n+j+1:2n)-\beta *A(:,n+j+1:n)w*w';$\\	
   	 		\indent \%Updating $S$ (if needed):\\
   	 		\indent\indent    $S(:,j+1:n)=S(:,j+1:n)-\beta *S(:,j+1:n)w*w';$\\	
   	 		\indent\indent    $S(:,n+j+1:2n)=S(:,n+j+1:2n)-\beta *S(:,n+j+1:n)w*w';$\\	
   	 		\indent \indent end	\\
   	\indent \indent end\\
  \indent 	end\\
   \end{alg}
   \begin{alg}\label{vlg}
   	function[c,s]=vlg(k,a)\\
   	\indent $twon=length(a);\;n=twon/2;$\\
   	\indent $r=\sqrt{a(k)^2+a(n+k)^2};$\\
   	\indent if $r=0$ then $c=1;\;s=0;$\\
   	\indent else
   	\indent $\displaystyle{c=\frac{a(k)}{r};\;\;s=\frac{a{(n+k)}}{r}};$\\
   	\indent end
   \end{alg}
   \begin{alg}\label{vlh}
   	function[$\beta$,w]=vlh(k,a)\\
   	\indent $twon=length(a);\;n=twon/2;$\\
   	\indent \% $w=(w_1,\ldots,w_{n-k+1})^T;$\\
   	\indent  $r1=\sum_{i=2}^{n-k+1} a(i+k-1)^2;$\\
   		\indent  $r=\sqrt{a(k)^2+r1};$\\
   	\indent $w_1= a(k)+sign(a(k))r;$\\
   	\indent $w_i=a{(i+k-1)}$ for $i=2,\ldots,n-k+1;$\\
   		\indent  $r={w_1^2+r1};\;\;
   	\displaystyle{\beta=\frac{2}{r};}$\\
   	\indent \%$ P=I-\beta ww^T;\;\;(H(k,w)a)_i=0$ for $i=k+1,\ldots,n.$\\
   	\indent end
   \end{alg}
\subsection{Numerical experiments}
In this work, we restrict our selves to the algorithmic aspect of  J Hessenberg reduction of a matrix, via symplectic Householder transformations. We showed how this reduction may be handled. The reduction process involves free parameters. We outlined how some optimal choice can be done, which gave rise to JHOSH
  algorithm. The latter uses only symplectic Householder transformations, which are not orthogonal. We succeed to replace half of them by transformations which are both orthogonal and symplectic. This gave rise to JHMSH algorithm, which behaves with satisfactory properties and is better than  all the  previous ones. 
   Very important questions on numerical aspects as for example the other choices of the free parameters,  breakdowns, near breakdowns, different strategies to cure these near breakdowns, and also their early  prediction before performing  computations which are not necessary, and so on, deserves a   detailed study. This will be the focus of a forthcoming paper. Nevertheless,
   we propose below two significant numerical examples in the following sense :   in the literature, to our knowledge,  only  the JHESS algorithm is used to perform a $J$-Hessenberg reduction of a matrix, with symplectic transformations. The JHESS belongs to the same class of algorithms as are JHOSH and JHMSH. The figures below compare JHMSH, JHMSH$2$ (which is a slight modification of JHMSH) and JHESS. The numerical examples show that the   later, as presented in \cite{Buns}   meets a fatal breakdown and  thus fails for all $n\geq 3$ , while  the  JHMSH, JHMSH$2$, with a slight modification,  work up with very satisfactory  precision.
   Let us consider the following matrix 
   $A=
   \begin{pmatrix}
   M_{11} & M_{12} \\ 
   M_{21}& M_{22}
   \end{pmatrix},
   $   
  with 
   $M_{11}=
   \begin{pmatrix}
   1 &   &   &   \\ 
   2 & 1 &   &   \\ 
   & \ddots  & \ddots  &   \\ 
   &   & 2 & 1
   \end{pmatrix} 
   $, 
   $M_{12} =
   \begin{pmatrix}
   1 & 2 &   &   \\ 
   2 & 1 & \ddots &   \\ 
   & \ddots & \ddots & 2 \\ 
   &   & 2 & 1
   \end{pmatrix} 
   $, $M_{21}=
   \begin{pmatrix}
   0 & 2 &   &   \\ 
   0 & 1 & \ddots &   \\ 
   &  \ddots & \ddots & 2 \\ 
   &   & 0 & 1
   \end{pmatrix} 
   $
    and 
   $M_{22} =
   \begin{pmatrix}
   1 &   &   &   \\ 
   3 & 1 &   &   \\ 
   & \ddots & \ddots &   \\ 
   &   & 3 & 1
   \end{pmatrix} . 
   $ Each block $M_{ij}$ is of size $n\times n.$ We obtain\\
   \begin{tabular}{|c|c|c|c||c|c|c|}
   	\hline
   	$n$ & \multicolumn{3}{|c}{Loss of $J$-Orthogonality $\left\Vert I-S^{J}S\right\Vert _{2}$} & 
   	\multicolumn{3}{||c|}{Error of the reduction $\left\Vert H-S^{J}AS\right\Vert _{2}$} \\ \hline
   	& $JHESS$ & $JHMSH$ & $JHMSH2$ & $JHESS$ & $JHMSH$ & $JHMSH2$ \\ \hline
   	$2$ & fails & $2.5168e-16$ & \multicolumn{1}{|c|}{$3.1402e-16$} & fails & $%
   	1.0361e-15$ & $1.0262e-15$ \\ \hline
   	$3$ & fails & $1.0412e-15$ & \multicolumn{1}{|c|}{$9.7146e-16$} & fails & $%
   	1.0623e-14$ & $5.6678e-15$ \\ \hline
   	$4$ & fails & $3.1015e-15$ & \multicolumn{1}{|c|}{$3.6572e-15$} & fails & $%
   	6.3153e-14$ & $2.9172e-14$ \\ \hline
   	$5$ & fails & $2.8250e-14$ & \multicolumn{1}{|c|}{$3.3284e-14$} & fails & $%
   	1.4279e-13$ & $6.8545e-14$ \\ \hline
   	$6$ & fails & $4.1918e-14$ & \multicolumn{1}{|c|}{$4.3812e-14$} & fails & $%
   	2.5845e-13$ & $1.6997e-13$ \\ \hline
   	$7$ & fails & $2.0709e-13$ & \multicolumn{1}{|c|}{$1.1965e-13$} & fails & $%
   	2.7021e-12$ & $5.7755e-13$ \\ \hline
   	$8$ & fails & $1.7497e-12$ & \multicolumn{1}{|c|}{$7.4477e-13$} & fails & $%
   	1.0972e-11$ & $3.5435e-12$ \\ \hline
   	$9$ & fails & $1.2988e-10$ & \multicolumn{1}{|c|}{$5.8035e-11$} & fails & $%
   	1.0461e-09$ & $3.8219e-10$ \\ \hline
   	$10$ & fails & $4.8062e-10$ & \multicolumn{1}{|c|}{$1.1476e-10$} & fails & $%
   	3.4164e-09$ & $7.1532e-10$ \\ \hline
   	$11$ & fails & $6.6942e-10$ & \multicolumn{1}{|c|}{$1.7784e-10$} & fails & $%
   	4.7274e-09$ & $5.7041e-10$ \\ \hline
   	$12$ & fails & $4.5165e-10$ & \multicolumn{1}{|c|}{$1.7250e-10$} & fails & $%
   	1.1306e-08$ & $8.0399e-10$ \\ \hline
   	$13$ & fails & $7.9908e-10$ & \multicolumn{1}{|c|}{$2.9785e-10$} & fails & $%
   	7.4063e-09$ & $1.7637e-09$ \\ \hline
   	$14$ & fails & $7.6406e-10$ & \multicolumn{1}{|c|}{$1.7497e-10$} & fails & $%
   	8.3607e-09$ & $1.0158e-09$ \\ \hline
   	$15$ & fails & $1.7248e-09$ & \multicolumn{1}{|c|}{$1.9073e-10$} & fails & $%
   	1.1932e-08$ & $9.8201e-10$ \\ \hline
   	$16$ & fails & $6.9530e-10$ & \multicolumn{1}{|c|}{$1.9133e-10$} & fails & $%
   	5.6770e-09$ & $1.1922e-09$ \\ \hline
   	$17$ & fails & $1.9515e-09$ & \multicolumn{1}{|c|}{$2.1889e-10$} & fails & $%
   	1.4054e-08$ & $1.2598e-09$ \\ \hline
   	$18$ & fails & $1.1824e-09$ & \multicolumn{1}{|c|}{$6.2781e-10$} & fails & $%
   	1.4967e-07$ & $5.7161e-09$ \\ \hline
   	$19$ & fails & $3.6906e-09$ & \multicolumn{1}{|c|}{$2.2293e-10$} & fails & $%
   	2.5400e-08$ & $1.4194e-09$ \\ \hline
   	$20$ & fails & $2.8172e-09$ & \multicolumn{1}{|c|}{$2.6019e-10$} & fails & $%
   	1.2725e-07$ & $2.0413e-09$ \\ \hline
   	$21$ & fails & $1.5606e-08$ & \multicolumn{1}{|c|}{$8.6765e-10$} & fails & $%
   	2.6936e-07$ & $5.1208e-09$ \\ \hline
   	$22$ & fails & $1.0522e-09$ & \multicolumn{1}{|c|}{$2.4081e-10$} & fails & $%
   	1.1047e-08$ & $1.9222e-09$ \\ \hline
   	$23$ & fails & $3.8242e-09$ & \multicolumn{1}{|c|}{$2.6805e-10$} & fails & $%
   	2.1954e-08$ & $1.6025e-09$ \\ \hline
   	$24$ & fails & $1.1119e-09$ & \multicolumn{1}{|c|}{$4.8392e-10$} & fails & $%
   	5.6800e-08$ & $3.2751e-09$ \\ \hline
   	$25$ & fails & $3.9755e-09$ & \multicolumn{1}{|c|}{$4.2710e-10$} & fails & $%
   	2.2816e-08$ & $2.6839e-09$ \\ \hline
   	$26$ & fails & $1.8132e-09$ & \multicolumn{1}{|c|}{$1.4496e-09$} & fails & $%
   	3.2416e-08$ & $1.0678e-08$ \\ \hline
   	$27$ & fails & $1.2417e-08$ & \multicolumn{1}{|c|}{$1.1257e-09$} & fails & $%
   	1.0768e-07$ & $1.0010e-08$ \\ \hline
   	$28$ & fails & $2.2564e-09$ & \multicolumn{1}{|c|}{$1.1255e-09$} & fails & $%
   	1.4462e-07$ & $8.2262e-09$ \\ \hline
   	$29$ & fails & $3.9904e-08$ & \multicolumn{1}{|c|}{$2.3791e-09$} & fails & $%
   	6.3257e-07$ & $4.1958e-08$ \\ \hline
   	$30$ & fails & $1.6554e-09$ & \multicolumn{1}{|c|}{$5.4776e-10$} & fails & $%
   	5.9380e-08$ & $4.0406e-09$ \\ \hline
   \end{tabular}
   \\
   \\
   \\
Consider now the Hamiltonian case : \\$A=\left(
\begin{array}{ll}
 M_{11} & M_{12} \\ 
 M_{21} & M_{22}
 \end{array}
\right),
$   
 where  
 $M_{11} = 
 \begin{pmatrix}
 1 &  &  &  \\ 
 2 & 1 &  &  \\ 
 & \ddots & \ddots &  \\ 
 &  & 2 & 1%
 \end{pmatrix}
 $, $M_{12} = 
 \begin{pmatrix}
 1 & 2 &  &  \\ 
 2 & 1 & \ddots &  \\ 
 & \ddots & \ddots & 2 \\ 
 &  & 2 & 1%
 \end{pmatrix}
 $,
 $M_{21} = 
 \begin{pmatrix}
 0 & 0 &  &  \\ 
 0 & 1 & 3 &  \\ 
 & 3 & \ddots & 3 \\ 
 &  & 3 & 1%
 \end{pmatrix}
 $  and  $M_{22} = -M_{11}^T.
 $
 We get \\
 \\
 \begin{tabular}{|c|c|c|c||c|c|c|}
 	\hline
 	$n$ & \multicolumn{3}{|c}{Loss of $J$-Orthogonality $\left\Vert I-S^{J}S\right\Vert _{2}$} & 
 	\multicolumn{3}{||c|}{Error of the reduction $\left\Vert H-S^{J}AS\right\Vert _{2}$} \\ 
 	\hline
 	& $JHESS$ & $JHMSH$ & $JHMSH2$ & $JHESS$ & $JHMSH$ & $JHMSH2$ \\ 
 	\hline
 	$2$ & fails & $1.3843e-16$ & $2.7756e-17$ & fails & $3.4732e-16$ & $%
 	7.5047e-16$ \\ \hline
 	$3$ & fails & $2.1967e-15$ & $4.1153e-15$ & fails & $1.4123e-14$ & $%
 	9.5826e-15$ \\ \hline
 	$4$ & fails & $3.1724e-14$ & $1.1623e-14$ & fails & $1.0235e-13$ & $%
 	1.1283e-13$ \\ \hline
 	$5$ & fails & $5.5639e-13$ & $4.5393e-13$ & fails & $2.2678e-12$ & $%
 	1.4082e-12$ \\ \hline
 	$6$ & fails & $1.3229e-14$ & $3.1824e-14$ & fails & $1.6308e-13$ & $%
 	1.8500e-13$ \\ \hline
 	$7$ & fails & $1.9456e-13$ & $2.9018e-13$ & fails & $4.2300e-12$ & $%
 	5.7276e-12$ \\ \hline
 	$8$ & fails & $2.4182e-13$ & $9.1255e-14$ & fails & $2.6360e-12$ & $%
 	1.2184e-12$ \\ \hline
 	$9$ & fails & $7.0030e-12$ & $4.6008e-12$ & fails & $2.8308e-11$ & $%
 	6.0019e-11$ \\ \hline
 	$10$ & fails & $6.7908e-11$ & $1.8421e-11$ & fails & $1.8128e-10$ & $%
 	4.2484e-11$ \\ \hline
 	$11$ & fails & $1.2746e-10$ & $3.6111e-11$ & fails & $1.2132e-09$ & $%
 	1.3393e-10$ \\ \hline
 	$12$ & fails & $1.6379e-09$ & $1.1448e-10$ & fails & $5.6804e-09$ & $%
 	1.0683e-09$ \\ \hline
 	$13$ & fails & $5.7401e-09$ & $1.8386e-09$ & fails & $4.3477e-07$ & $%
 	5.7596e-09$ \\ \hline
 	$14$ & fails & $5.9220e-09$ & $2.7826e-09$ & fails & $1.1117e-07$ & $%
 	1.1405e-08$ \\ \hline
 	$15$ & fails & $1.1198e-07$ & $1.5282e-08$ & fails & $8.4815e-07$ & $%
 	2.1596e-07$ \\ \hline
 	$16$ & fails & $3.2853e-07$ & $1.9260e-07$ & fails & $3.6979e-06$ & $%
 	8.2332e-07$ \\ \hline
 	$17$ & fails & $1.0707e-06$ & $1.9526e-07$ & fails & $1.4713e-05$ & $%
 	3.9805e-06$ \\ \hline
 	$18$ & fails & $2.2014e-04$ & $2.2887e-05$ & fails & $1.3000e-03$ & $%
 	4.6621e-04$ \\ \hline
 	$19$ & fails & $7.0710e-05$ & $2.0118e-05$ & fails & $1.5000e-03$ & $%
 	4.0607e-04$ \\ \hline
 	$20$ & fails & $7.9995e-04$ & $4.0086e-05$ & fails & $4.1000e-03$ & $%
 	6.8321e-04$ \\ \hline
 	 \end{tabular}
\\
\\
\section{Conclusion}
In this paper, we presented a reduction of a matrix to the upper $J$-Hessenberg form, based on the symplectic Householder transformations, which are rank-one modification of the Identity. This reduction is the crucial step for constructing an efficient SR-algorithm. The method is the analog of the reduction of a matrix to Hessenberg form, via Householder transformations, when instead of an Euclidean linear space, one takes a sympletctic one.  Then the algorithm JHOSH is derived, corresponding to an optimal choice of the free parameters. Furthermore, JHOSH is significantly improved by showing that half of these symplectic Householder transformations may be replaced by Van Loan's symplectic and orthogonal transformations leading to two variants JHMSH and JHMSH2 which are significantly more stable numerically. The numerical experiments confirm  the expected results.


\end{document}